\documentclass[11pt,oneside]{amsart}
\usepackage{epsfig}
\usepackage{amsmath,ifthen, amsfonts, amssymb, srcltx, amsopn}
\usepackage{graphicx}
\usepackage[small]{caption}

\newcommand{\showcomments}{yes}
\renewcommand{\showcomments}{no}

\newsavebox{\commentbox}
\newenvironment{com}%
{\ifthenelse{\equal{\showcomments}{yes}}%
{\footnotemark
    \begin{lrbox}{\commentbox}
    \begin{minipage}[t]{1.25in}\raggedright\sffamily\tiny
    \footnotemark[\arabic{footnote}]}
{\begin{lrbox}{\commentbox}}}%
{\ifthenelse{\equal{\showcomments}{yes}}%
{\end{minipage}\end{lrbox}\marginpar{\usebox{\commentbox}}}
{\end{lrbox}}}

\newtheorem{thm}{Theorem}[section]
\newtheorem{lem}[thm]{Lemma}
\newtheorem{mlem}[thm]{Main Lemma}
\newtheorem{cor}[thm]{Corollary}

\newtheorem*{main}{Main Theorem}

\theoremstyle{definition}
\newtheorem{defn}[thm]{Definition}

\newtheorem{rem}[thm]{Remark}
\newtheorem{exmp}[thm]{Example}

\newtheorem{conv}[thm]{Convention}

\newcommand{\isomto}{\overset{\simeq}{\rightarrow}}

\newcommand{\field}[1]{\mathbb{#1}}
\newcommand{\integers}{\ensuremath{\field{Z}}}

\newcommand{\Euclidean}{\ensuremath{\field{E}}}

\newcommand{\boundary}   {{\ensuremath \partial}}

\begin{document}

\title[C(6) groups do not contain $F_2 \times F_2$]
{C(6) groups do not contain $F_2 \times F_2$}

\author[H.~Bigdely]{Hadi Bigdely}
      \address{Dept. of Math \& Stats.\\
               McGill University \\
               Montreal, Quebec, Canada H3A 2K6 }
\email{bigdely@math.mcgill.ca}
\author[D. T.~Wise]{Daniel T. Wise}
      \email{wise@math.mcgill.ca}
\thanks{Research supported by NSERC}


\date{\today}

\begin{com}
{\bf \normalsize COMMENTS\\}
ARE\\
SHOWING!\\
\end{com}

\begin{abstract}
We show that a group with a presentation satisfying the C(6) small cancellation condition cannot contain a subgroup isomorphic to $F_2 \times F_2$.
\end{abstract}

\maketitle

\tiny
\normalsize
\section{Introduction}
The goal of this paper is to prove the following result:
\begin{main}
Let $G$ be a group with a presentation satisfying the C(6) small-cancellation condition.
Then $G$ does not contain a subgroup isomorphic to  $F_2 \times F_2$.
\end{main}

We recall that a 2-complex $X$ satisfies the C($p$)-T($q$) small-cancellation condition if each reduced disc diagram $D\rightarrow  X$ has the property that its internal 2-cells of $D$ have at least~$p$  bordering $2$-cells (locally) and that internal 0-cells of $D$ have (either $2$ or) at least $q$ adjacent $2$-cells. We note that the $T(3)$ conditions holds for any 2-complex (and so we often write C($p$)
instead of C($p$)-T(3)).
A presentation satisfies the C($p$)-T($q$) condition if its standard 2-complex does.
In a certain sense, the C($p$)-T($q$) condition represents a ``combinatorial comparison'' condition with a simply-connected surface
tiled by $p$-gons with $q$~meeting around each vertex. The cases of greatest interest are when $\frac1p+\frac1q\leq \frac12$, when the corresponding tiling corresponds to a regular tiling of the Euclidean  or Hyperbolic plane.

For a finite 2-complex $X$, the group  $\pi_1X$ is readily seen to be word-hyperbolic if it satisfies C($p$)-T($q$) with
$\frac1p+\frac1q<\frac12$
\ (i.e. C(7)-T(3),  C(5)-T(4),  C(4)-T(5), \&  C(3)-T(7)).
However, $\pi_1X$ only necessarily manifests features of nonpositive curvature when  $\frac1p+\frac1q=\frac12$
\ (i.e. C(6)-T(3), C(4)-T(4), \& C(3)-T(6)).
For instance, Gersten and Short showed that $\pi_1X$ is automatic in this case \cite{GerstenShort91}.

It is well-known that a word-hyperbolic group cannot contain a $\integers\times \integers$ subgroup,
and in this setting, the $\integers\times \integers$ subgroup in $\pi_1X$ leads to a combinatorial flat plane in $\widetilde X$.
More generally, failure of word-hyperbolicity corresponds to failure of a linear isoperimetric function which corresponds to the existence of a combinatorial flat plane in $\widetilde X$  as shown by Ivanov and Schupp \cite{IvanovSchupp98}.
However the degree to which $\pi_1X$ fails to be word-hyperbolic has not yet been studied deeply.

One sense in which $\pi_1X$ can ``strongly fail'' to be hyperbolic is if there is a profusion of $\integers\times \integers$ subgroups,
or indeed, if these subgroups richly ``interact'' with each other, as in $F_2\times F_2$.
This can certainly occur when $X$ is C(4)-T(4), as indeed $F_2\times F_2 \cong \pi_1X$ when $X= B\times B$ where $B$ is a bouquet of 2~circles.
A conclusion of this work is that C(6)-T(3) groups appear closer to being hyperbolic than C(4)-T(4) groups.
We do  not report on this here in detail, but the class of C(3)-T(6) groups also cannot contain $F_2\times F_2$ for similar reasons.
\begin{com}This will be explained in detail and within a broader framework in \cite{BigdelyThesis}.
\end{com}

We now give a brief description of the sections of the paper.
In Section~\ref{sec:sc} we review the definitions of small-cancellation theory that we will need.
In Section~\ref{sec:l.c} we examine locally convex maps and the properties of a certain thickening of a locally convex subcomplex.
In Section~\ref{sec:bitorus} we study locally convex maps $Y\rightarrow X$ that are associated with conjugacy classes
of $F_2\times \integers$ in $\pi_1X$. These are the main objects of interest in the paper.
In Section~\ref{sec:mr} we prove the main result.


\section{Small-cancellation theory}\label{sec:sc}
In this section we give a brief review of the basic notions of
small-cancellation theory. We follow the geometric language given in \cite{McCammondWiseFanLadder},
and more  details and examples can be
found there. A more classical reference is \cite{LS77}.

We shall work in the category of combinatorial complexes:
\begin{defn}[Combinatorial maps and complexes]\label{defn:combinatorial}
A map $Y\rightarrow X$ between CW complexes is \emph{combinatorial} if
its restriction to each open cell of $Y$ is a homeomorphism onto an
open cell of $X$.  A CW complex $X$ is \emph{combinatorial}
if the attaching map of each open cell of $X$ is combinatorial for a suitable subdivision.
\end{defn}

\begin{defn}[Disc diagram]\label{defn:diagrams}
A {\em disc diagram} $D$ is a compact contractible $2$-complex
with a fixed embedding in the plane. A {\em boundary cycle} $P$
of $D$ is a closed path in $\partial D$ which travels entirely
around $D$ (in a manner respecting the planar embedding of $D$).

A \emph{disc diagram in $X$} is a map $D\rightarrow X$.
It is a well-known fact, due to van Kampen, that whenever $P\rightarrow X$ is a nullhomotopic closed path,
there is a disc diagram $D\rightarrow X$ such that $P\rightarrow X$ factors as $P\rightarrow \boundary D\rightarrow X$.

Let $R_1$ and $R_2$ be $2$-cells that meet along a $1$-cell $e$ in the disc diagram $D\rightarrow X$. We say $R_1$ and $R_2$ are a \emph{cancellable pair} if the boundary paths of $R_1$ and $R_2$ starting at $e$ map to the same closed path in $X$. $D\rightarrow X$ is \emph{reduced} if it has no cancellable pair of 2-cells.
\end{defn}

\begin{defn}[Piece]\label{defn:piece}
Let $X$ be a combinatorial $2$-complex. Intuitively, a piece of
$X$ is a path which is contained in the boundaries of the
$2$-cells of $X$ in at least two distinct ways.  More precisely, a
nontrivial path $P\rightarrow X$ is a \emph{piece} of $X$ if there
are $2$-cells $R_1$ and $R_2$ such that $P\rightarrow X$ factors
as $P \rightarrow R_1 \rightarrow X$ and as $P\rightarrow
R_2\rightarrow X$ but there does not exist a homeomorphism
$\partial R_1\rightarrow \partial R_2$ such that there is a
commutative diagram:
\begin{equation*}\label{eq:piece}
\begin{array}{ccc}
P             & \rightarrow & \partial R_2\\
\downarrow    & \nearrow    & \downarrow\\
\partial R_1 & \rightarrow & X
\end{array}
\end{equation*}
Excluding commutative diagrams of this form ensures that $P$
occurs in $\partial R_1$ and $\partial R_2$ in essentially
distinct ways.
\end{defn}

\begin{defn}[C($p$)-complex]\label{defn:sc}
An \emph{arc} in a disc diagram is a path whose internal vertices have
valence~$2$ and whose initial and terminal vertices have
valence~$\geq 3$. The arc is {\em internal} if its interior lies
in the interior of $D$, and it is a {\em boundary arc} if it lies
entirely in $\boundary D$. A $2$-complex $X$ satisfies the C($p$) condition if the boundary path of each
$2$-cell in each reduced disc diagram $D$ either contains a nontrivial boundary arc,
 or is the concatenation of at
least~$p$ nontrivial internal arcs. A group $G$ is C(6)
if it is the fundamental group of a C(6) 2-complex.
\end{defn}

\begin{figure}[h]
\includegraphics[width=1 in]{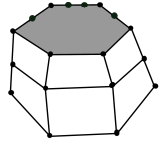}
\caption{
 The shaded $2$-cell~$R$ is
$3$-shell of $D$.} \label{fig:outercell}
\end{figure}

\begin{defn}[$i$-shell, spur]\label{defn:i-shells}
Let $D$ be a diagram.  An \emph{$i$-shell} of $D$ is a $2$-cell $R
\hookrightarrow D$ whose boundary cycle $\boundary R$ is the
concatenation $P_0P_1\cdots P_i$ where $P_0 \rightarrow D$ is a
boundary arc, the interior of $P_1\cdots P_i$ maps to the interior
of $D$, and  $P_j \rightarrow D$ is a nontrivial interior arc of
$D$ for all $j > 0$. The path $P_0$ is the {\em outer path} of the
$i$-shell. Note that $P_0 = \partial R \cap \partial D$. (See Figure~\ref{fig:outercell}.)

A $1$-cell $e$ in $\boundary D$ that is incident with a valence~$1$ $0$-cell $\upsilon$ is a {\em {spur}}.
\end{defn}

The following result
 is ``Greendlinger's Lemma'' which is the fundamental tool of small-cancellation theory
\cite[Thm~9.4]{McCammondWiseFanLadder} (see \cite[Thm~V.4.5]{LS77} for a classical version of this statement).
\begin{thm}[Greendlinger's Lemma]\label{green}
If $D$ is a C(6)-$T(3)$ disc diagram, then
 either $D$ is a single 0-cell or a single closed embedded 2-cell,
or else $D$ has at least $2\pi$ worth of spurs and $i$-shells with $i\leq 3$,
where each spur, 0-shell, and 1-shell contributes $\pi$,
each $2$-shell contributes  $\frac{2\pi}3$,  and each $3$-shell contributes~$\frac\pi3$.
\end{thm}

\begin{defn}[Missing $i$-shell]
Consider the commutative diagram on the left
$$\begin{array}{ccc}
P &\rightarrow & Y\\
\downarrow &  & \downarrow\\
R & \rightarrow & X\\
\end{array}
\hspace{2cm}
\begin{array}{ccc}
P &\rightarrow & Y\\
\downarrow &\nearrow  & \downarrow\\
R & \rightarrow & X\\
\end{array}
$$
where  $Y$ is a $2$-complex, $R$ is a closed $2$-cell, and
$P\rightarrow X$ is a path which factors through both $Y$ and $R$.
Let $\boundary R=PS$ where $S$ is the concatenation of $i$ pieces.
We say that $R$ is a {\em missing $i$-shell } for $Y$ if the map
$P\rightarrow Y$ does not extend to a map $R\rightarrow Y$ so that
the diagram on the right commutes.
\end{defn}

\begin{figure}[h]
\includegraphics[width=2.5 in]{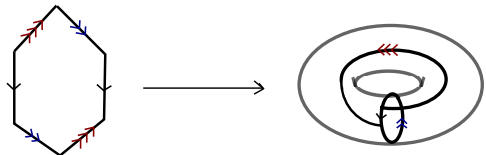}
\caption{Identifying opposite sides of a hexagon yields a hexagonal torus.} \label{fig:hextorus}
\end{figure}

\begin{defn}[Hexagonal torus and Honeycomb]
A {\em honeycomb in $X$} is a hexagonal tiling of $\Euclidean^2$ with some valence~2 vertices added.
A {\em hexagonal torus} is a C(6)-complex homeomorphic to a torus.
A very simple hexagonal torus is indicated in Figure~\ref{fig:hextorus}.
Of course, any hexagonal torus is the quotient of a honeycomb by a free cocompact action of $\integers \times \integers$.
\end{defn}

\section{Locally convex maps}\label{sec:l.c}

In this section we define locally convex and strongly locally convex maps and we show that the ``thickening'' of a strongly locally convex subcomplex is also strongly locally convex.
We  will follow the following convention for the remainder of the paper:
\begin{conv}\label{conv}
When we state $X$ is a $2$-complex, we mean a C(6)-complex. We
will assume that no two $2$-cells of $X$ have the same attaching
map.
\end{conv}

\begin{defn}[Locally convex map]
A combinatorial map
between C(6)-complexes is an {\em{immersion}} if it is locally
injective. An immersion $\phi:Y \rightarrow X$ between C(6)-complexes
is {\emph{locally convex}} if it does not have a missing $i$-shell
where $1\leq i\leq3$.
Figure~\ref{fig:locally convex} illustrates an example in which
the inclusion of the shaded subcomplex $Y$ is not a locally convex
map, indeed it has a missing $3$-shell $R$.

\begin{figure}[h]
\includegraphics[width=1 in]{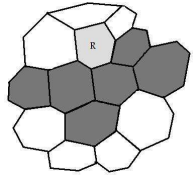}
\caption{$R$ is a $3$-shell in the complement of the dark shaded complex and is missing along its outer path in the shaded complex. So the inclusion map from the dark complex to the whole complex is not locally convex.} \label{fig:locally convex}
\end{figure}
\end{defn}

\begin{lem}
A locally convex map $f: Y\rightarrow X$ between simply connected
C(6)-complexes is injective.
\end{lem}

\begin{proof}
Let $\delta$ be a non closed path in $Y$ which maps to a closed
path $\gamma$ in $X$. Since $X$ is simply connected, $\gamma$
bounds a disc diagram. Choose $\delta$ among all pathes between
endpoints of $\delta$ such that its image $\gamma$ has the minimal
area disc diagram $D$ in $X$. By Lemma~\ref{green}, $D$ has an $i$-shell called $R$ where
$1\leqslant i\leqslant 3$. Let $\boundary R= QS$ where $S$ is the
concatenation of $i$ pieces ($1 \leqslant i \leqslant3$) and $Q$
is the outer path of the shell. Since $Y\rightarrow X$ is locally
convex, the map $R\rightarrow X$ induces a map $R\rightarrow Y$
otherwise, $R$ will be a missing $i$-shell for $Y$ . We show the
image of $R$ in $Y$ by $R$. So $Q$ is part of $\delta$. Now push
$Q$ toward $S$ in $Y$ to get a new path $\acute{\delta}$ whose end
points are the same as $\delta$. The image of $\acute{\delta}$ in
$X$ bounds the disc diagram $\acute{D}=D-R$ where
Area($\acute{D}$)$=$ Area($D$)$-1$ which is contradiction.
\end{proof}

\begin{lem}\label{intersection}
Let $X$ be a C(6)-complex. Let $Y_1$ and $Y_2$ be subcomplexes of $X$ such that each inclusion map $Y_i \hookrightarrow X$ is locally convex. Then $Y_1\cap Y_2 \hookrightarrow X$ is also locally convex.

\end{lem}
\begin{proof}
This follows immediately from the definition.
\end{proof}
\begin{defn}[Strongly locally convex subcomplex]
Let $X$ and $Y$ be C(6)-complexes. We call $\widetilde Y\hookrightarrow \widetilde X$ \emph{strongly locally convex} if for any $2$-cell $R$ with $\bar{R} \cap \widetilde Y\neq \varnothing$, either $R\subseteq \widetilde Y$ or each component of $\bar{R} \cap \widetilde Y$ is the concatenation of at most $2$ pieces. For example it is immediate that $\bar R$ is strongly locally convex whenever $R$ is a single $2$-cell. Honeycombs are also strongly locally convex.

Observe that strongly locally convex implies locally convex. Consequently, if $\widetilde Y \subseteq \widetilde X$ is strongly locally convex then $\bar{R}\cap \widetilde Y$ is actually connected by Lemma~\ref{intersection}.

We emphasize that the definition requires that components of $\bar{R}\cap \widetilde Y$ be expressible as the concatenation of at most $2$ pieces. It is possible that they are also expressible as the concatenation of more than two pieces.
\end{defn}

\begin{lem}{\label{trisect}}
If $\widetilde Y \subseteq \widetilde X$ is strongly locally convex and $R_1$, $R_2$ are $2$-cells in $\widetilde X$ with $\bar R_1 \cap \bar R_2 \neq \varnothing$, $\bar R_1 \cap \widetilde Y \neq \varnothing$ and $\bar R_2 \cap \widetilde Y \neq \varnothing$ then $\bar R_1 \cap \bar R_2 \cap \widetilde Y \neq \varnothing$.
\end{lem}

\begin{proof}
We show that if $\bar R_1 \cap \bar R_2 \neq \varnothing$, $\bar R_1 \cap \widetilde Y \neq \varnothing$ and $\bar R_2 \cap \widetilde Y \neq \varnothing$ then $\bar R_1 \cap \bar R_2$ is a singleton or a piece that intersects  $\widetilde Y$. Observe that $\bar R_1 \cap \bar R_2$ has one component. Assume $\bar R_1 \cap \bar R_2$ does not intersect $\widetilde Y$. Let $D$ be the minimal area disc diagram whose boundary path consists of the pathes $\alpha$, $\beta$ and $\gamma$ where $\alpha \subseteq \boundary \bar R_1$, $\beta \subseteq \boundary \bar R_2$ and $\gamma \subseteq \widetilde Y$. The disc diagram $D$ is illustrated as the dark complex in Figure~\ref{fig:twocells}-A.

\begin{figure}[h]
\includegraphics[width=3.3 in]{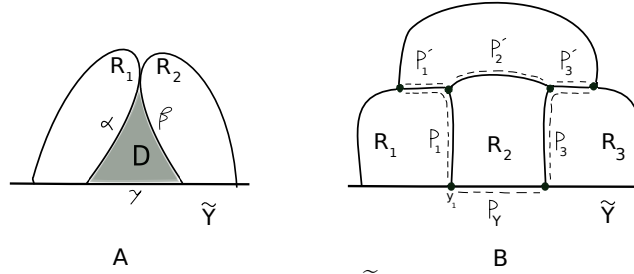}
\caption{In Figure~A, the complex $\widetilde Y$ is strongly locally convex and $\bar {R_1}, \bar {R_2}$ and $\widetilde Y$ should triply intersect. In Figure~B, the $2$-cells $R_1, R_2, R_3$ are subsets of ${\sf N}(\widetilde Y)$ and $\bar {R} \cap {\sf N}(\widetilde Y)$ is the concatenation of $P_1' P_2' P_3'$. }\label{fig:twocells}
\end{figure}
Observe that since the boundary path of $D$ has at most $4$~pieces, $D$ is not a single $2$-cell. By Lemma~\ref{green}, since $D$ does not have spurs, it must contain at least three $i$-shells where $i\leq 3$. Also since $\widetilde Y$ is strongly locally convex, if $D$ has $i$-shells, they must lie in the corners. $D$ can not have an $i$-shell in the corner corresponding to $\bar R_1$ and $\bar R_2$, therefore $D$ has at most two $i$-shells which is contradiction.
\end{proof}

We now show that a ``nice extension'' of a strongly locally convex subcomplex of a C(6)-complex, is again strongly locally convex.

\begin{defn}[Thickening]
Let $X$ be a C(6)-complex and $\widetilde Y \subseteq \widetilde X$ be a subcomplex. The \emph{thickening ${\sf N} (\widetilde Y)$ of $\widetilde Y$} is the subcomplex $$
{\sf N}(\widetilde Y)= \widetilde Y \cup \{\bar{R}~|~R~  \mbox {is a $2$-cell and}~ \bar{R}\cap \widetilde Y\neq  \varnothing\}.$$ We use the notation ${\sf N}^0 (\widetilde Y)= \widetilde Y$ and ${\sf N}^{i+1} (\widetilde Y)={\sf N}({\sf N}^{i}(\widetilde Y))$. Note that ${\sf N}(\widetilde Y)={\sf N}^1(\widetilde Y)$ might not contain an open neighborhood of $\widetilde Y$.
\end{defn}

\begin{rem}\label{union}
If $\widetilde X$ is connected and has no isolated $1$-cell and $\widetilde Y \neq \varnothing$ then $\widetilde X=\cup_{k\geq 0} {\sf N}^k(\widetilde Y)$. Indeed for any path $P$ whose initial vertex is on a cell $\alpha$ in $\widetilde X$ and whose terminal vertex lies on $\widetilde Y$, we see that $\alpha \subset {\sf N}^{|P|}(\widetilde Y)$.
\end{rem}

\begin{lem}\label{thickening}
Let $X$ be a C(6)-complex and $\widetilde Y \subseteq \widetilde X$ be a connected subcomplex. If $\widetilde Y\hookrightarrow \widetilde X$ is strongly locally convex then
${\sf N} (\widetilde Y) \hookrightarrow \widetilde X$ is also strongly locally convex.

\end{lem}
\begin{proof}
Let $R$ be a $2$-cell in $\widetilde X$ such that $\bar R \cap {\sf N}(\widetilde Y) \neq \varnothing$. Suppose a subpath $P$ of $\bar R \cap {\sf N}(\widetilde Y)$ is the concatenation $P_1' P_2' P_3'$ where each $P_i'$ is a path in $\bar R\cap \bar R_i$ and each $R_i\subseteq {\sf N}(\widetilde Y)-\widetilde Y$. The lemma follows easily from the following claim:
The subpath $P$ can be expressed as the concatenation of at most two pieces.

Proof of the claim: Without loss of generality assume that $\bar R_1 \cap \bar R_2$ and $\bar R_1 \cap \bar R_3$ are both nonempty. By Lemma~\ref{trisect} $\bar R_1$, $\bar R_2$ and $\widetilde Y$ triply intersect, also $\bar R_3$, $\bar R_2$ and $\widetilde Y$ triply intersect. Let $P_1$ be the shortest path containing $\bar R_1 \cap \bar R_2$ from a point $y_1$ in $\bar R_1 \cap \bar R_2 \cap \widetilde Y$ to the initial point of $P_1' P_2' P_3'$. Similarly let $P_3$ be the shortest path containing $\bar R_3 \cap \bar R_2$ from the terminal point of $P_1' P_2' P_3'$to a point $y_3$ in $\bar R_2 \cap \bar R_3 \cap \widetilde Y$. Let $P_Y$ be a path in $\bar R_2 \cap \widetilde Y$ between $y_3$ and $y_1$. Consider the path $P_1(P_1' P_2' P_3')P_3P_Y= (P_1P_1') (P_2')(P_3'P_3)P_Y$. By hypothesis $P_Y$ is at most two pieces and thus (after removing the backtracks in $P_1P_1'$ and $P_3'P_3$) the path in $\boundary \bar R_2$ is the concatenation of less than $6$ pieces. Therefore the path can not travel around $R_2$ and thus travels through an arc $A$ in $\boundary \bar R_2$. We claim that $R \subset {\sf N}(\widetilde Y)$ otherwise $P_1'P_2'P_3'\cap P_Y=\varnothing$ and therefore $P_1$ and $P_3$ must intersect in $A$. But then $\bar R_1 \cap \bar R_3 \neq \varnothing$ and by Lemma~\ref{trisect}, $\bar R_1 \cap \bar R_3 \neq \varnothing$ which implies that $P_1'P_2'P_3'$ is replaceable by $P_1''P_3''$. We refer the reader to Figure~\ref{fig:twocells}-B.
\end{proof}

\section{Bitorus}\label{sec:bitorus}
In this section, we define $2$-complexes called ``bitori'' which are the main objects of interest in the paper. Also, we study locally convex maps $Y\rightarrow X$ that are associated with conjugacy classes
of $F_2\times \integers$ in $\pi_1X$.

\begin{defn}[Band and Slope]
Let $X$ be a honeycomb in which all pieces have length $1$. Two
edges are {\em{equivalent}} if they are antipodal edges of a $2$-cell
in $X$. This generates an equivalence relation for $1$-cells of $X$. A \emph{band} is a minimal subcomplex of $X$
containing an equivalence class. Note that a band corresponds to a sequence of hexagons inside a honeycomb where
attaching $1$-cells are antipodal.
In a honeycomb we have three families of bands. Each
band has two boundaries which we call {\em{slopes}}.
So we have three different families of slopes. Let $X$ be a complex whose universal cover $\widetilde X$ is a flat plane. An {\em{immersed band}} in $X$ is the image of a band by the covering map. Note that interior of a band in $X$ embeds but it is possible for slopes to get identified. Also note that bands do not cross themselves. Two distinct slopes are parallel if they do not cross.
\end{defn}

A \emph{flat annulus} is a concentric union of $n\geq 0$ bands. Equivalently, it is the complex obtained from a hexagonal torus by removing a single band.

\begin{figure}[h]
\includegraphics[width=3.5 in]{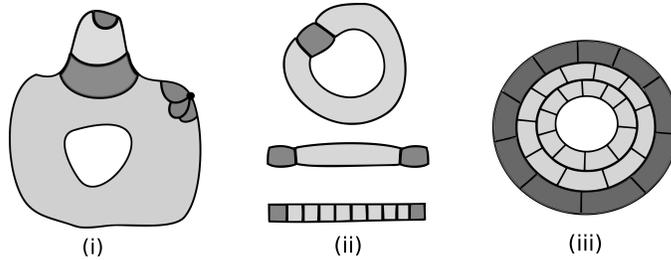}
\caption{Figure~(i) illustrates an $i$-shell which arises if a 2-cell is multiply external, and a (negatively curved) valence~4 vertex on $\boundary A$. Figure~(ii) illustrates the conclusion that can be drawn if a 2-cell contains an edge in both boundary paths.
 Figure~(iii) illustrates the outer band of $4$-shells that is sliced off to obtain a smaller annular diagram.} \label{fig:annulus}
\end{figure}

The following can be proven along the same lines as proofs in \cite{McCammondWiseFanLadder}:
\begin{lem}\label{flat annulus}
Let $A$ be a compact nonsingular annular $C(6)$ $2$-complex with no spurs or $i$-shells with $i\leq 3$ along either of its boundary paths.
Then $A$ is a flat annulus.
\end{lem}
\begin{proof}
We assign a $\frac{2\pi}{3}$ angle to each internal corner of valence~$\geq 3$, a $\frac{\pi}{2}$ angle to each corner with a single boundary edge, and a $\pi$ angle to all other corners.
All internal $0$-cells and all internal $2$-cells have curvature $\leq 0$.
No closed $2$-cell $R$ intersects the same boundary path of $A$ in two or more disjoint subpaths,
since by Theorem~\ref{green} there would then be an $i$-shell with $i\leq 3$ in $A$ at a subdiagram of $A$ subtended by
$R$, as indicated in Figure~\ref{fig:annulus}-(i).
If some $2$-cell $R$ intersects both boundary paths of $A$, then by cutting along $R$, we obtain a disc diagram $L$ with at most two
$i$-shells (with $i\leq 3$) and hence $L$ is a nonsingular ``ladder'' by \cite[Thm~9.4]{McCammondWiseFanLadder},
and consequently $A$ was a ``one band annulus'' to begin with. See Figure~\ref{fig:annulus}-(ii).

Now we show that all $0$-cells and $2$-cells have curvature exactly~$0$. By Theorem~\cite[Thm~4.6]{McCammondWiseFanLadder}, we have:
\begin{align}\label{combinatorial GB}
\sum_{f\in \text{2-cells(A)}} \text{curvature}(f)+\sum_{v\in \text{0-cells(A)}} \text{curvature}(v)&=2\pi \cdot \chi(A)
\end{align}
If there is an $i$-shell with $i\geq 5$ in one of the boundary paths, or there is an interior $2$-cell with more than $6$ pieces or
 an internal or external $0$-cell of valence~$\geq 4$ then the left side of  Equation~ \eqref{combinatorial GB} would be negative, but this would contradict that the right side is~$0$.  
  Consider the 2-cells whose boundaries contain an edge in the outside boundary path of $A$.
  Since each of these 2-cells forms a $4$-shell and since there are no valence~$4$ vertices on $\boundary A$,
  we see that consecutive such 2-cells meet in a nontrivial piece, and this sequence of 2-cells forms a width~1 annular band.
  The subdiagram obtained by removing this band is again a nonsingular annular diagram with no $i$-shell with $i\leq 3$,
  for otherwise $A$ would have had an internal 2-cell with $\leq 5$ sides.
  The result now follows by induction, and $A$ is a  union of bands as in Figure~\ref{fig:annulus}-(iii).
\end{proof}

\begin{defn}[bitorus]
A \emph{bitorus} is a compact and connected $C(6)$-complex homeomorphic with $B \times S^{1}$ where $B$ is a finite connected leafless graph and $\chi (B)=-1$. There are three families of these complexes:
The first family which is homeomorphic to a complex constructed by attaching a flat annulus to two tori along some slope, is the union of bands
attached along parallel slopes. (As mentioned there are three families of slopes in a torus). Figure~\ref{fig:bitorus}-(i) illustrates an example of this family where the attaching slope does not wrap around the torus. The second family is homeomorphic to a complex constructed by attaching two tori along some slope. Figure~\ref{fig:bitorus}-(ii) illustrates this family but in general, a slope can wrap around a torus several times. The third is homeomorphic to the $2$-complex obtained by attaching a flat annulus to a torus along two parallel and separate slopes.
\begin{figure}[h]
\includegraphics[width=3.8 in]{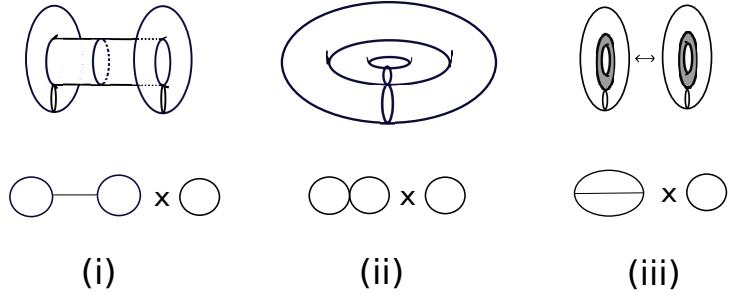}
\caption{We illustrate the three types of bitori. Figure~(iii) is a bitorus obtained by attaching two hexagonal tori along the shaded regions which is the union of bands.
 } \label{fig:bitorus}
\end{figure}

\end{defn}

We now prove a lemma that plays an important role in the main
theorem:

\begin{mlem}\label{convex in bitorus}
Let $X$ be a bitorus. Let $Y$ be a compact and connected C(6)-complex and $f : Y\rightarrow X$ a combinatorial map which is
$\pi_1$-injective and locally convex. Then either $\pi_1Y\cong 1$ or $\pi_1Y\cong \integers$ or $
\integers \times \integers \subseteq \pi_1Y$.
\end{mlem}
\begin{exmp}
The statement of Lemma~\ref{convex in bitorus} does not hold if we replace the ``bitorus'' by an analogous complex $Z$ that is constructed from three tori instead of two such that $\pi_1Z \cong \langle a_1,a_2,a_3,a_4~|~[a_i,a_{i+1}]=1 : 1\leq i \leq4 \rangle$. Indeed let $Z$ be the $2$-complex obtained by attaching the $2$-cells $A$, $B$ and $C$ to $Z^1$, as indicated in Figure~\ref{fig:exm}. Let $Y$ be the graph indicated in Figure~\ref{fig:exm} and observe that the inclusion map $i:Y\hookrightarrow Z$ is locally convex.
\begin{figure}[h]
\includegraphics[width=5.3 in]{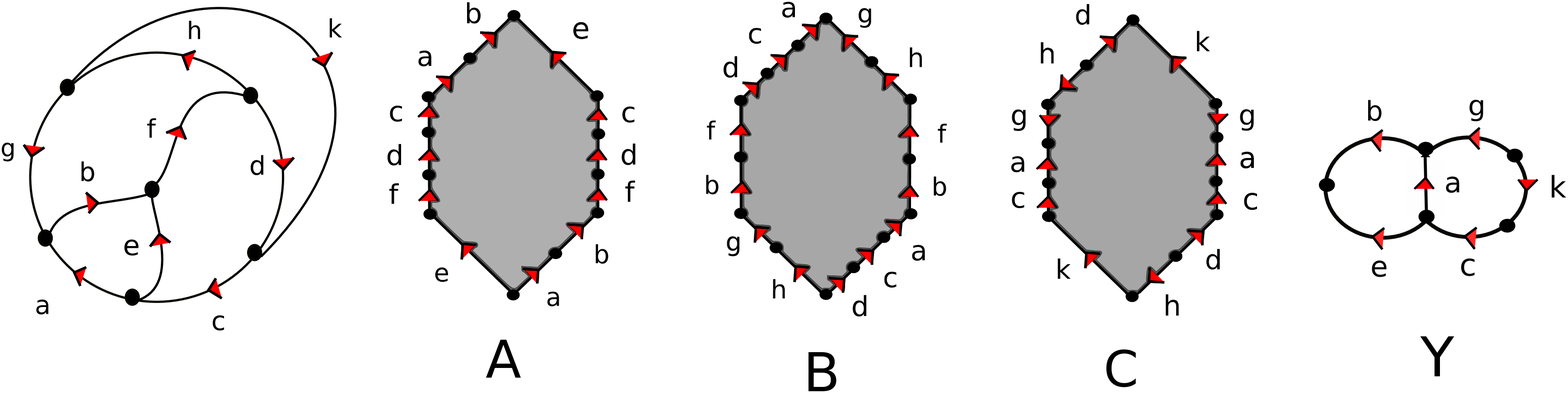}
\caption{After attaching the $2$-cells $A$, $B$ and $C$ to $Z^1$ we get a C(6)-complex $Z$ and the inclusion map $i:Y\hookrightarrow Z$ is locally convex.} \label{fig:exm}
\end{figure}
\end{exmp}
\begin{exmp}
For $C(4)$-$T(4)$ complexes, an immersion is {\em{locally convex}} if it has no missing $i$-shell for $i=0,1,2$. But Lemma~\ref{convex in bitorus} fails in this case. Indeed $F_2 \times F_2 \cong \pi_1 X$ where $X$ is the $C(4)$-$T(4)$ complex that is the product of two graphs.
\end{exmp}

\begin{proof}[Proof of Lemma~\ref{convex in bitorus}]

Consideration of all circles that are in the same slope as attaching circles yields
a graph of spaces $\Gamma_{X}$ whose vertex spaces are circles and
 whose edge spaces are bands. Figure~\ref{fig:cylinder}-$B$ illustrates $\Gamma_{X}$. Let $X_v$ and $X_e$ be respectively a vertex space and an edge space where $v\in \Gamma^0_{X}$ and $e\in \Gamma^1_{X}$. The graph of spaces for $X$ will induce a
graph of spaces $\Gamma_{Y}$ for $Y$ where:

$$Y_v=f^{-1} (X_v)~~~~~~~~~\hspace{0.3 in} ~~~~~~~~~ \hspace{0.3 in} Y_e=f^{-1} (X_e)$$

First, assume there is a vertex space in $\Gamma_{Y}$ which is a circle
called $C$. If there is no $2$-cell attaching to $C$, since $f$ is
$\pi_1$-injective and locally convex, $Y=C$ and $\pi_1Y=\integers$. Otherwise since there is no missing $3$-shell, the edge space attached to $C$ is a
band and therefore all edge spaces in $Y$ are bands. Figure~\ref{fig:cylinder}-$A$ illustrates an example of cylindrical edge space. In this case, if $\Gamma_{Y}$ contains a circle then $
\integers \times \integers \subseteq \pi_1Y$ and if $\Gamma_{Y}$ does not contain a circle then $Y$ is homotopy equivalent to a circle and $\pi_1Y\cong \integers$.
\begin{figure}[h]
\includegraphics[width=2.1 in]{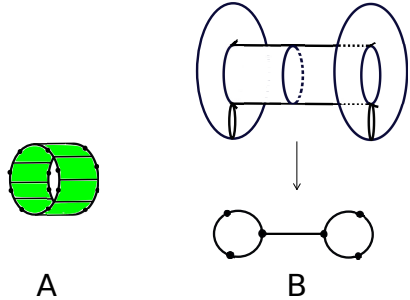}
\caption{Figure~A illustrates an edge space which is a band and in Figure~B, the $2$-complex is obtained by attaching a flat annulus to two tori and each torus is a union of three bands. } \label{fig:cylinder}
\end{figure}

We have a map $g : Y\rightarrow \Gamma_{Y}$.

Case 1: $Y$ does not contain a $2$-cell and no vertex space of $\Gamma_{Y}$ is a circle. Specifically each vertex space is each vertex space is a point or a subcomplex of a circle which is not closed. We show that $\Gamma_{Y}$ has no valence $3$ vertex. If there exists a valence $3$ vertex in
$\Gamma_Y$, then the image of $Y$ by $f$ locally looks like the
dark path in $1$-skeleton of $X$ in Figure~\ref{fig:valence3}. In
this case we will have a missing $3$-shell which contains $0$-cell
with valence~$3$ in $X$ and this is contradiction. So in this
case, the valence of each vertex in $\Gamma_Y$ is $\leq 2$. Either $\Gamma_Y$ is a circle and $\pi_1Y \cong \integers$ or $\Gamma_Y$ is not a circle in which
 $Y$ is contractible and $\pi_1Y \cong 1$.
\begin{figure}[h]
\includegraphics[width=1.8 in]{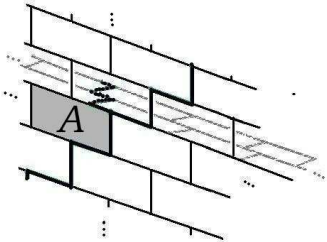}
\caption{A is a missing $3$-shell for the dark complex inside $X$.} \label{fig:valence3}
\end{figure}

Case 2: $Y$ contains some $2$-cells and no vertex space of $\Gamma_{Y}$ is a circle. We show that
$\pi_1Y \cong \integers$ or $1$. First note that since there is no
missing $3$-shell, the difference between the number of $2$-cells on
two adjacent edge spaces of $\Gamma_Y$ is at
most one. Assume a length~$3$ path in the graph $\Gamma_Y$ then the corresponding $2$-cells of edge spaces can not retreat and then extend. Therefore in a length $3$ path in $\Gamma_Y$ the number of corresponding $2$-cells of edge spaces can not decrease and then increase. Figure~\ref{fig:inc-dec}-$B$ illustrates an example in which the number of the edge spaces corresponding to a path in $\Gamma_Y$ retreat and then extend. If there is
no valence $3$ vertex in $\Gamma_Y$ then either $\Gamma_Y$ is a circle and $\pi_1Y\cong \integers$ or $\Gamma_Y$ does not contain a circle and
$\pi_1Y \cong\integers$.

\begin{figure}[h]
\includegraphics[width=2.8 in]{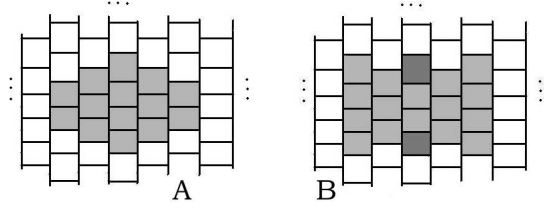}
\caption{In Figure~A, the number of $2$-cells in edge spaces
decreases after increasing (from right to left), so there is no
missing $3$-shell. But, in Figure~B, this number first decreases
and then increases and we have two very dark missing $3$-shells.}
\label{fig:inc-dec}
\end{figure}

Now, we will discuss the case in which we have a valence~$3$
vertex in $\Gamma_Y$. Consider $3$ edge spaces $\varepsilon_1$, $\varepsilon_2$ and $\varepsilon_3$ meeting at a vertex space $\nu$. Assume that they have the same number of
$2$-cells $n$. We know that the vertex space $\nu$ is a segment of
a circle. Therefore it has two vertices of valence $1$ called $\omega_1$ and $\omega_2$. One of the edge spaces $\varepsilon_1$, $\varepsilon_2$ and $\varepsilon_3$ contains none of $\omega_1$ and $\omega_2$. Assume $\varepsilon_1$ does not contain $\omega_1$ and $\omega_2$. In Figure~\ref{fig:edge-number}-$A$ the dark edge in the back is $\varepsilon_1$.

\begin{figure}[h]
\includegraphics[width=3.4 in]{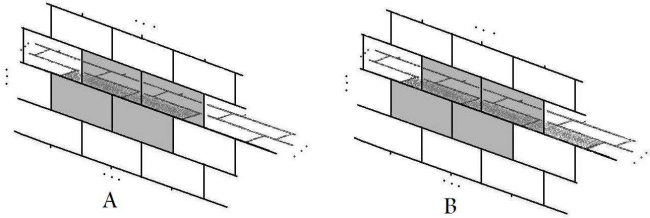}
\caption{Three edge spaces $\varepsilon_1, \varepsilon_2, \varepsilon_3$ meet along a vertex space $\nu$.
The rear edge space is $\varepsilon_1$ Figure~A, and $\varepsilon_3$ in Figure~B.}\label{fig:edge-number}
\end{figure}

Considering $\varepsilon_1$ and $\varepsilon_2$, the edge space $\varepsilon_1$ is retreating in one side, also considering $\varepsilon_1$ and $\varepsilon_3$, the edge space $\varepsilon_1$ is retreating. Therefore we should not have missing $3$-shell, the
number of $2$-cells in the edge space attached to $\varepsilon$
which does not have intersection with $\nu$ is $n-1$. Also,
since the number of $2$-cells of edge spaces in this branch
decreased from $n$ to $n-1$, it should decrease in the next edge spaces. So
we showed that if we have a valence $3$ vertex in $\Gamma_Y$
called $\upsilon$ and all edge spaces attached to $\upsilon$ have
the same number of $2$-cells in $Y$, then in one of the paths
(branches) called $\tau$ when we travel far from $\upsilon$, the
number of $2$-cells of edge spaces will decrease.
 Therefore in
this case, the image of $\tau$ in $\Gamma_Y$ will not result a loop in $\Gamma_Y$ and we can collapse $\tau$
without changing $\pi_1Y$.

Now, assume the $3$ edge spaces $\varepsilon_1$, $\varepsilon_2$ and $\varepsilon_3$ having intersection in the
vertex space $\nu$, do not have the same number of $2$-cells. So
two of them should have the same number of
$2$-cells $m$ and the third one $m-1$ or $m+1$. Let $\varepsilon_1$ and $\varepsilon_2$ have $m$ $2$-cells. If $\varepsilon_3$
has $m-1$ $2$-cells, then by the same argument in the previous case, we can collapse this branch without changing the fundamental group. Now assume that $\varepsilon_3$ has
$m+1$  $2$-cells. Since the number of $2$-cells from the
$\varepsilon_3$ to $\varepsilon_1$ and $\varepsilon_2$ decreases, the image
of both of them in $\Gamma_Y$ are not part of a
loop. Figure~\ref{fig:edge-number}-$B$ illustrates this case.

Therefore if case~1 or case~2 occurs then either $\Gamma_Y$ is homotopy equivalent to a circle in which case $\pi_1Y \cong \integers$, or
$Y$ is contractible so $\pi_1Y \cong 1$.
\end{proof}
We will  employ Lemma~\ref{convex in bitorus} in the following contrapositive form:

\begin{cor}\label{free not convex}
There is no locally convex, $\pi_1$-injective map $Y\rightarrow X$ where $X$ is a bitorus and $Y$ is a compact connected $2$-complex with $\pi_1Y \cong F_2$.
\end{cor}
\begin{lem}\label{immersion}
Let $f:\widetilde Y\rightarrow \widetilde X$ be a map where $Y$ is a bitorus and $X$ is a C(6)-complex. Then $f$ is strongly locally convex. In particular $f$ is locally convex.
\end{lem}
\begin{proof}
We first show that $f$ is locally convex.
Suppose that $Q\rightarrow \widetilde Y$ is the outer path of a missing $i$-shell $R$ with $i\leq 3$.
Since the inner path of $R$ is the concatenation of $i\leq 3$ pieces, the C(6) condition applied to $R$ shows that  $Q$ cannot be the concatenation of $\leq 2$~pieces in $\widetilde Y$. It follows that $Q$ must fully contain two consecutive maximal pieces in
 the boundary of a single $2$-cell $R'$ of $\widetilde Y$. But this violates the C(6) condition for $R'$.

 Having proven local convexity, we turn to strong local convexity.
Suppose $R$ is a $2$-cell that is not in $\widetilde Y$ such that $\bar R\cap \widetilde Y \neq \varnothing$.
Assume that $P=\bar R\cap \widetilde Y$ cannot be expressed as the concatenation of at most two pieces.
As before, consideration of paths in the honeycomb $\widetilde Y$, shows that $P$ contains two consecutive maximal pieces in the boundary of a single 2-cell $R'$. This violates the C(6) condition for $R'$.
 \end{proof}

\begin{lem}\label{constructin l.c}
Let $X$ be a C(6)-complex such that $F_2 \times \integers \subseteq \pi_{1}X$. There exists a $2$-complex $Y'$ equal to $Y\vee Q$ where $Y$ is a bitorus, $Q=[0,n]$ and $0$ is identified with a $0$-cell in $Y$ and $n$ is the basepoint. And there exists a basepoint preserving map $f:Y'\rightarrow X$ such that $f \mid_ {Y}$ is locally convex and the following diagram commutes:
\begin{equation*}
\begin{array}{ccc}
              &             & \pi_1 Y'\\
              & \reflectbox{\rotatebox[origin=c]{315}{$\isomto$}}    & \downarrow\\
F_2\times \integers & \hookrightarrow & \pi_1 X
\end{array}
\end{equation*}
\end{lem}

\begin{proof}
 We will construct $Y'=Y \vee Q$ and an immersion $f:Y'\rightarrow X$, by Lemma~\ref{immersion}, $f \mid_{Y}$ is locally convex.

 Let $\nu \in X^0$ be the basepoint. Let $F_2 \times \integers \cong \langle a_1,a_2,c~|~[a_1,c] , [a_2,c] \rangle.$ For $i=1,2$, let $A_i$ and $C_i$ be closed, based paths in $X$, such that $A_i$ represents $a_i$ and $C_1$ and $C_2$ both represent $c$. For $i=1,2$, let $D_i$ be a minimal area disc diagram with boundary path $A_iC_iA_i^{-1}C_i^{-1}$. Moreover we shall make the above choices such that $D_i$ has minimal area among all such choices of $A_i$ and $C_i$. By identifying the top and bottom $C_i$ paths and identifying the left and right $A_i$ paths we obtain a quotient $T_i$ of $D_i$. Observe that $T_i={T_i}' \vee [0,n_i]$ is the wedge of a torus with $[0,n_i]$. Moreover there exists an induced combinatorial map $T_i\rightarrow X$. The minimality of $D_i$ ensures that ${T_i}' \rightarrow X$ is an immersion. Let $V_i=[0,n_i]$. By possibly folding, we can shorten $V_i$ to assume that $T_i \rightarrow X$ is also an immersion. Note that our original paths $C_i\rightarrow X$ corresponds to two paths $C_i \rightarrow D_i$ which are then identified to a single path $C_i \rightarrow T_i$ which we shall now examen. For $i=1, 2$, let $U_i$ be an embedded closed path in ${T_i}'$ such that  $C_i=V_iU_iV_i^{-1}$. The complex $T_i$ is illustrated in Figure~\ref{fig:torustail}-A.

\begin{figure}
\includegraphics[width=3 in]{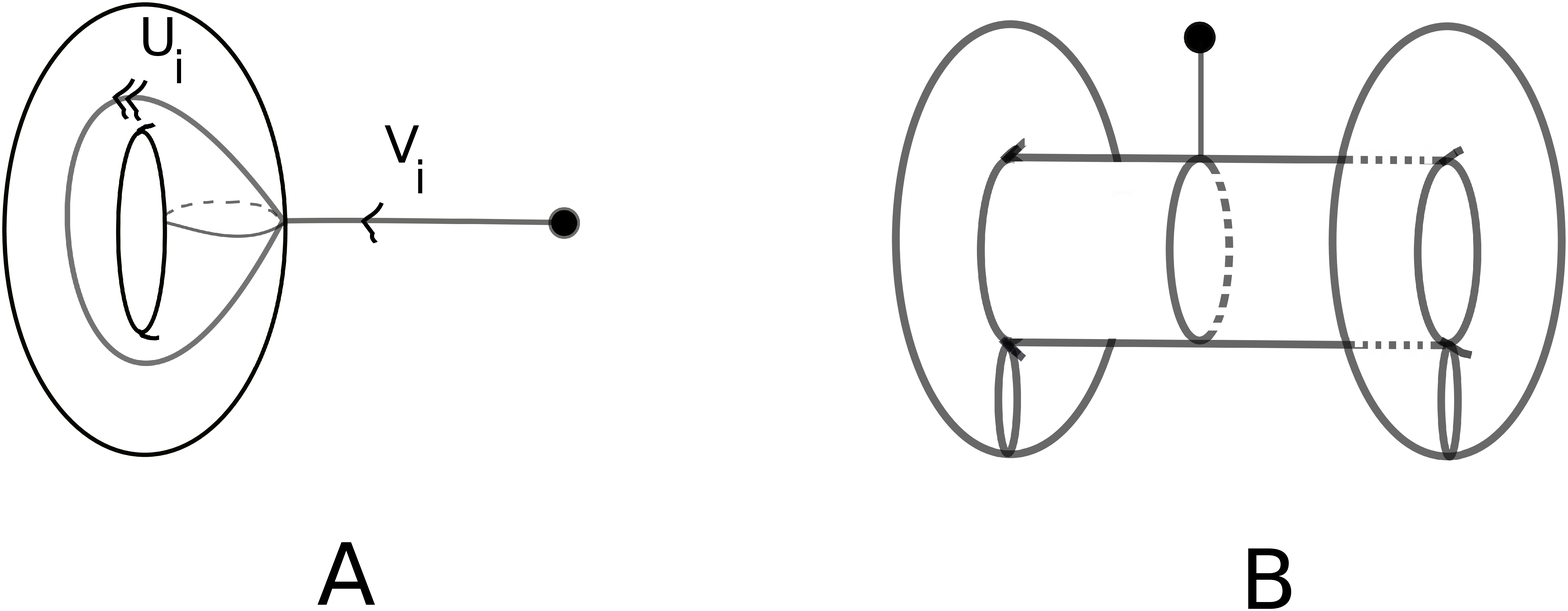}
\caption{On the left is $T_i$ which is an hexagonal torus with an arc attached to it. On the right is $N$.}\label{fig:torustail}
\end{figure}

Let $A$ be an annular diagram whose boundary paths $P_1$, $P_2$ are respectively homotopic to the image of $U_1$ in $U_1 \rightarrow {T_1}' \rightarrow X$ and the image of $U_2$ in $U_2 \rightarrow {T_2}' \rightarrow X$ and whose conjugator is path-homotopic to the image of ${V_1}^{-1}V_2$ in $V_1^{-1}V_2 \rightarrow X$. Moreover choose $A$ such that it has minimal area with these properties. Note that $A$ is non-singular since $U_1$ is path-homotopic to $U_2$.
Consider the base lifts of $\widetilde{T_1}$ and $\widetilde{T_2}$ to $\widetilde{X}$. Note that these determine lifts of $\widetilde{{T_1}'}$ and $\widetilde{{T_2}'}$. Either $\widetilde{{T_1}'}$ and $\widetilde{{T_2}'}$ intersect or do not intersect.

 We first consider the case where $\widetilde{{T_1}'}$ and $\widetilde{{T_2}'}$ do not intersect. Observe that for $0 \leq i \leq 3$, $A$ has no missing $i$-shell along either boundary path. Indeed ${T_1}'$ and ${T_2}'$ do not have missing $i$-shell and so if $A$ had an $i$-shell, we could reduce its area. By Lemma~\ref{flat annulus},  $A$ is a flat annulus. Let $Y$ be the $2$-complex obtained by attaching $A$ to ${T_1}' \sqcup {T_2}'$ along $P_1, P_2$. Note that the annulus is nonsingular and it is a strip whose one side is identified in ${T_1}'$ and the other is identified in ${T_2}'$. Since $A$ has minimal area, by the above argument, there is no folding and $S$ is a bitorus of the first type. Moreover since the conjugator of $A$ is path-homotopic to ${V_1}^{-1}V_2$, there is a path $Q=[0, n]$ where $0$ is identified in $A$ and $n$ is identified with the basepoint $n_1=n_2$ in $V_1$ and $V_2$. Figure~\ref{fig:torustail}-B shows the $2$-complex $Y \vee Q$. In conclusion, in this case $Y'$ equals  $Y \vee Q$ where $Y$ is a bitorus of the first type and $Q$ corresponds to a basepath $[0,n]$ where $0$ is identified in $A$ and $n$ is the basepoint.

We now consider the case where $\widetilde{{T_1}'}$ and $\widetilde{{T_2}'}$ intersect. By Lemma \ref{intersection}, $\widetilde{{T_1}'} \cap \widetilde{{T_2}'}$ is a locally convex subcomplex of $\widetilde X$ and it is also infinite since the element $c$ stabilizes both $\widetilde{{T_1}'}$ and $\widetilde{{T_2}'}$. Observe that since $\widetilde{{T_1}'} \cap \widetilde{{T_2}'}$ is locally convex and infinite, it is a slope or union of consecutive bands. Therefore $\widetilde{{T_1}'} \cap \widetilde{{T_2}'}$ contains a periodic line. Let
$A_1=\widetilde A/ c$ and let $S={T_1}' \sqcup {T_2}' / A_1$.

If there is no folding between ${T_1}'$ and ${T_2}'$, then $A_1$ is a slope and $Y'$ equals $Y \vee Q$ where $Y=S$ is a bitorus of the second type and $Q$ corresponds to a basepath $[0,n]$ where $0$ is identified in the common slope and $n$ is the basepoint. Otherwise, we start to fold ${T_1}'$ with ${T_2}'$. Note that the folding process cannot identify ${T_1}'$ with ${T_2}'$, since $\pi_1{T_1}'$ and $\pi_1{T_2}'$ are not commensurable in $\pi_1X$. By considering all slopes parallel to a given slope, it is natural to regard ${T_i}'$ as a graph of spaces whose vertex spaces are circles and whose edge spaces are bands. If a $2$-cell $R_1$ in ${T_1}'$ is folded  with a $2$-cell $R_2$ in ${T_2}'$ then the entire band containing $R_1$ will be folded  with the band containing $R_2$. Moreover, note that the C(6) structure of $S$ ensures that these bands consist of the same number of $2$-cells. As a result ${T_1}'$ and ${T_2}'$ will be identified along a union of consecutive bands and we call the obtained complex $Y$. In conclusion the result of folding process in this case is a complex $Y'$ equal to $Y \vee Q$ where $Y$ is a bitorus of the third type and $Q$ corresponds to a basepath $[0,n]$ where $0$ is identified in one of ${T_1}'$ or ${T_2}'$ and $n$ is the basepoint. Moreover in all cases, there exists an induced combinatorial map $f: Y'\rightarrow X$ such that $f \mid_{Y}$ is an immersion
and therefore locally convex by Lemma~\ref{immersion}.
\end{proof}

\section{Main Result}\label{sec:mr}


\begin{thm}\label{main them1}
A C(6) group cannot contain $F_2 \times F_2$.
\end{thm}
\begin{proof}
Let $X$ be a based C(6)-complex whose fundamental group is $G$. Suppose that $G$ contains $F_2
\times F_2 \cong \langle a,b \rangle \times \langle c,d \rangle$, we will reach a contradiction.
Without loss of generality, we can assume that each 1-cell of $X$ lies on a 2-cell.
Indeed, since we are arguing by contradiction, we can replace $X$ by a smallest $\pi_1$-injective subcomplex
$X_o$ whose fundamental group contains $\langle a,b \rangle \times \langle c,d \rangle$. Since $F_2 \times F_2$ does not split as a free product, if $X_o$ contained a 1-cell not on the boundary of a 2-cell, then we could pass to a smaller $\pi_1$-injective subcomplex
whose fundamental group contains $F_2\times F_2$.

Consider the subgroups $G_1=\langle a,b \rangle \times \langle c\rangle \cong F_2 \times \integers$ and
 $G_2=\langle a,b \rangle \times \langle d \rangle \cong F_2
\times \integers$. Since $G_i \subset G$, by Lemma~\ref{constructin l.c}, there exists a $2$-complex ${Y_i}'$ that equals $Y_i\vee Q_i$ where $Y_i$ is a bitorus, $Q_i=[0,n_i]$ and $n_i$ is identified with a point in $Y_i$ and $0$ is the basepoint. Moreover, there exists a  basepoint preserving map $f:{Y_i}'\rightarrow X$ such that $f \mid_{Y_i}$ is locally convex and the following diagram commutes:
\begin{equation*}
\begin{array}{ccc}
              &             & \pi_1 {Y_i}'\\
              &\reflectbox{\rotatebox[origin=c]{315}{$\isomto$}}   & \downarrow\\
G_i & \hookrightarrow & \pi_1 X
\end{array}
\end{equation*}
By Remark~\ref{union}, there exists $k$ such that $\widetilde{{Y_1}} \cap {\sf N}^k(\widetilde{{Y_2}}) \neq \varnothing$.

Since by Lemma~\ref{thickening}, $\widetilde {{Y_1}} \rightarrow \widetilde X$ and  ${\sf N}^k(\widetilde{{Y_2}})  \rightarrow \widetilde X$ are both  locally convex, by Lemma~\ref{intersection}, $\widetilde {{Y_1}} \cap {\sf N}^k(\widetilde{{Y_2}})  \rightarrow \widetilde X$
is also locally convex.
Now observe that $\langle a,b \rangle=stab(\widetilde {{Y_1}}) \cap stab( {\sf N}^k(\widetilde{{Y_2}}) ) \subseteq stab(\widetilde {{Y_1}} \cap {\sf N}^k(\widetilde{{Y_2}})) $. Moreover the quotient space $Z=(stab(\widetilde {{Y_1}} \cap {\sf N}^k(\widetilde{{Y_2}})) \setminus (\widetilde {{Y_1}} \cap{\sf N}^k(\widetilde{{Y_2}}))) $ is compact, since $Z$ is a component of the fibre product of the maps ${Y_1}\rightarrow X$ and ${{\sf N}^k({Y_2})}\rightarrow X$. Since $F_2 \cong \langle a,b\rangle \subseteq \pi_1 Z$ and $Z \rightarrow {Y_i}'$ is locally convex with ${Y_i}'$ a bitorus, this contradicts Corollary~\ref{free not convex}.
\end{proof}

\bibliographystyle{alpha}
\bibliography{wise}

\def\cprime{$'$} \def\polhk#1{\setbox0=\hbox{#1}{\ooalign{\hidewidth
  \lower1.5ex\hbox{`}\hidewidth\crcr\unhbox0}}} \def\cprime{$'$}
  \def\cprime{$'$} \def\polhk#1{\setbox0=\hbox{#1}{\ooalign{\hidewidth
  \lower1.5ex\hbox{`}\hidewidth\crcr\unhbox0}}}
\begin{thebibliography}{MW02}

\bibitem[GS91]{GerstenShort91}
S.~M. Gersten and H.~Short.
\newblock Small cancellation theory and automatic groups. {I}{I}.
\newblock {\em Invent. Math.}, 105(3):641--662, 1991.

\bibitem[IS98]{IvanovSchupp98}
S.~V. Ivanov and P.~E. Schupp.
\newblock On the hyperbolicity of small cancellation groups and one-relator
  groups.
\newblock {\em Trans. Amer. Math. Soc.}, 350(5):1851--1894, 1998.

\bibitem[LS77]{LS77}
Roger~C. Lyndon and Paul~E. Schupp.
\newblock {\em Combinatorial group theory}.
\newblock Springer-Verlag, Berlin, 1977.
\newblock Ergebnisse der Mathematik und ihrer Grenz\-gebiete, Band 89.

\bibitem[MW02]{McCammondWiseFanLadder}
Jonathan~P. McCammond and Daniel~T. Wise.
\newblock Fans and ladders in small cancellation theory.
\newblock {\em Proc. London Math. Soc. (3)}, 84(3):599--644, 2002.

\end{thebibliography}
\end{document}